\documentclass[12pt]{amsart}
\usepackage{graphicx,amsmath}
\usepackage{xcolor}
\newtheorem{thm}{Theorem}[section]
\newtheorem{cor}[thm]{Corollary}
\newtheorem{lem}[thm]{Lemma}

\theoremstyle{definition}
\newtheorem{defn}[thm]{Definition}
\theoremstyle{remark}

\numberwithin{equation}{section}

\def\dis{\displaystyle}
\def\bq{\begin{equation}}
\def\eq{\end{equation}}
\def\bp{\begin{proof}}
\def\ep{\end{proof}}
\def\noi{\noindent}
\def\bi{\bigskip}
\def\me{\medskip}
\def\sm{\smallskip}
\def\P{{\mathbb P}^n}

\def\S{{ S^{k}_{ij}}}

\def\D{{\mathbb D}}
\def\al{\alpha}
\def\va{\varphi}
\def\ga{\gamma}
\def\si{\sigma}

\def\C{{\mathbb C}}
\def\B{{\mathbb B}^n}
\def\Cn{{\mathbb C}^n}
\def\la{{\lambda}}
\def\laa{{\lambda}_{\alpha}}
\def\mar{{\mu_r(\alpha)}}
\def\Fa{{\mathcal F}_{\alpha}}
\def\Far{{\mathcal F}_{\alpha, r}}
\def\Fao{{\mathcal F}_{\alpha}^0}
\def\I{\textsl{I}}

\begin{document}

\title{Families of Homomorphic Mappings in the Polydisk}

\author{Martin Chuaqui and Rodrigo Hern\'andez}
\thanks{The authors were partially supported by Fondecyt Grants 1150115,
1150284, 1190380, and
\endgraf
1190756.
\endgraf
\endgraf
{\sl MCS 2020}. Primary: 32H02, 32A17;\, Secondary: 30C45.}

\begin{abstract}
We study classes of locally biholomorphic mappings defined in the polydisk $\P$ that have bounded Schwarzian operator in the Bergman metric. We establish important properties of specific solutions of the associated system of differential equations, and show a geometric connection between the order of the classes and a covering property. We show for modified and slightly larger classes that the order is Lipschitz continuous with respect to the bound on the Schwarzian, and use this to estimate the order of the original classes.

\end{abstract}
\maketitle
\section{Introduction}
In this paper, we investigate families of locally biholomorphic mappings in the polydisk $\P$ defined by the Schwarzian operator. This is a classical theme in one complex variable, where the simply-connected case corresponds to mapping in the unit disk $\D$ with bounded Schwarzian in the hyperbolic norm. The nature of the group of automorphisms in $\D$ yields families that are {\sl linearly invariant}, with results that are often invariant under the group. The original papers by Pommerenke have become a landmark in this \cite{Po}. In several variables, linear invariance is present for the families in the ball $\B$, and allows for a fairly detailed analysis that includes a generalization of the Ahlfors-Weill extension and examples of locally linearly convex domains \cite{CH3},\cite{APS, Ho}. Our aim here is to study the case for the polydisk, a domain for which the explicit nature of the automorphisms and the Bergman metric partially compensates for the absence of linear invariance. Nevertheless, we were unable reproduce a {\sl density function} associated to the mappings treated that was both invariant and convex relative to the Bergman metric, missing thus a crucial step for the aforementioned Ahlfors-Weill extension.

In Section 2 we present preliminary facts about the Schwarzian operator and specific bounds that are satisfied in $\P$. In Section 3 we present the important notion of {\it normalization} that gives a geometric insight for the norm of the families. We also develop techniques based on comparison to study specific solutions of the linear system related to the Schwarzian operator. In particular, for sufficiently small values of the parameter in the family, we show that certain normalized solutions at the origin do vanish in $\P$. This can be thought as having the Ahlfors-Weill extension at one point, and implies that for $f$ in the family, the image $f(\P)$ misses the complex hyperplane $1-(a_1w_1+\cdots+a_nw_n)=0$, where $(n+1)(a_1,\ldots,a_n)=(\nabla J_f)(0)$, $J_f$ the complex jacobian of $f$. In one complex variable, this corresponds to the fact that $-2/f''(0)\notin f(\D)$ when $(1-|z|^2)^2|Sf(z)|\leq 2$.

In the final section we estimate the {\sl order} of the families treated. In \cite{CH2}, a generalization of Pommerenke's variational methods was employed to estimate the order of the families in $\B$, but linear invariance was required. Here instead, we appeal to a dilation argument that yields Lipschitz continuity for the order of slightly larger families. Some of our methods can be considered extensions of techniques that in one variable yield optimal results \cite{CO, EK, Ne}, and although sharp constants in $\Cn$ are in many cases still far out of reach, we have made an effort to be as precise as possible. The comparison techniques we develop are likely applicable on more general domains, and our results fit into a broader class of mappings in several variables studied by different authors, see for example \cite{HK, Ki, Pf, Su}.

\section{Preliminaries}

Let $f:\Omega\subset\Cn\rightarrow\Cn$ be a locally biholomorphic mapping
defined on some domain $\Omega$. There are several equivalent definitions of a Schwarzian derivative in several complex variables \cite{MM1, MM2, MT}, and we will follow the notion offered in \cite{Od} because of its explicit connection with differential equations.
The author defines a family of Schwarzian derivatives of $f=(f_1,\ldots,f_n)$ as
$$
\S f= \dis\sum^{n}_{l=1}\frac{\partial^{2} f_{l}}{\partial z_{i}\partial
z_{j}} \frac{\partial z_{k}}{\partial
f_{l}}-\frac{1}{n+1}\left(\delta^{k}_{i}\frac{\partial}{\partial
z_{j}}+\delta^{k}_{j}\frac{\partial}{\partial z_{i}}\right)
\log(J_f)\, ,
$$
where $i,j,k=1,2, \ldots,n,$
and $\delta^{k}_{i}$ are the Kronecker symbols. For $n>1$ the Schwarzian derivatives
have the properties $\S f=0$ for
all $i,j,k=1,2,\ldots ,n$ iff $f(z)= M(z),$ for some M\"{o}bius transformation
$$
M(z)=\left(\frac{l_1(z)}{l_0(z)},\ldots,\frac{l_n(z)}{l_0(z)}\right)\, ,
$$
where $l_i(z)=a_{i0}+a_{i1}z_1+\cdots+a_{in}z_n$ with
$\det(a_{ij})\neq 0$. Furthermore, for a composition one has
\bq\label{chain rule}
\S(g\circ f)(z)=\S f(z)+ \sum^{n}_{l,m,r=1}S^{r}_{lm}g(w)\frac{\partial w_{l}}{\partial z_{i}}
\frac{\partial w_{m}} {\partial z_{j}} \frac{\partial z_{k}} {\partial w_{r}}\; , \; w=f(z) \, .
\eq
Thus, if $g$ is a M\"{o}bius transformation then $\S(g\circ f)= \S f$. The coefficients $S^0_{ij}f$
are given by
$$
S^0_{ij}f(z)=(J_f)^{\frac{1}{n+1}}\left(\frac{\partial^2}{\partial z_i\partial z_j}(J_f)^{-\frac{1}{n+1}}-
\sum_{k=1}^n\,\frac{\partial}{\partial z_k}(J_f)^{-\frac{1}{n+1}}S^k_{ij}f(z)\right)\, .
$$
In \cite{Od} we can also find a description of the functions with prescribed Schwarzian derivatives $\S f$.
Consider the following over-determined system of partial differential equations,
\bq\label{system}
\frac{\partial^{2}u}{\partial z_{i}\partial z_{j}} =
\sum^{n}_{k=1}P^{k}_{ij}(z)\frac{\partial u}{\partial z_{k}}+
P^{0}_{ij}(z)u \; , \quad i,j= 1,2, \ldots,n\, ,
\eq
where $z=(z_{1},z_{2},...,z_n)\in \Omega$ and $P^k_{ij}(z)$ are
holomorphic functions in $\Omega$, for $i,j,k=0,\ldots,n$. The system (\ref{system})  is called
{\it completely integrable} if there are\,  $n+1$ (maximun)
linearly independent solutions, and is said to be in {\it canonical
form} (see \cite{Yo}) if the coefficients satisfy
$$
\sum_{j=1}^{n} P_{ij}^{j}(z)=0 \; ,\quad i=1,2,\ldots,n\, .
$$
It was shown that (\ref{system}) is a completely integrable system in canonical form
if and only if $P^k_{ij}=\S f$ for a locally biholomorphic mapping $f=(f_1, \ldots, f_n)$,
where $f_i=u_i/u_0$ for $1\leq i\leq n$ and $u_0, u_1,\ldots, u_n$ is a set of linearly independent
solutions of the system. For a given mapping $f$, $u_0=(J_f)^{-\frac{1}{n+1}}$ is always
a solution of (\ref{system}) with $P_{ij}^k=\S f$.

We recall the following definitions from \cite{He}, where the individual Schwarzians $S^kf$ are conformed to define the following operators.
\begin{defn}
For each $k=1,\ldots,n$ we let $S^{k}_f$ be the $n\times n$ matrix $$ S^k_f= (\S f)\, ,\quad i,j=1,\ldots,n\, .$$
\end{defn}

\begin{defn}
Let $f:\Omega\rightarrow\Cn$ be locally holomorphic, and let $T_z\Omega$ stand for the tangent space at $z\in\Omega$. We define the \textit{Schwarzian derivative operator} as the mapping $S_f(z):T_z\Omega \to T_z\Omega$ given by
$$
S_f(z)(\vec{v})=\left(S^1_f(z)(\vec{v})\, ,\ldots,S^{n}_f(z)\vec{v}\right)\, ,
$$ where
$\vec{v}\in T_z \Omega$ and $ S^k_f(z)(\vec v)=\vec{v}^{\,t}S^k_f(z)\vec{v}$.
\end{defn}
\noi With this, we can rewrite (\ref{system}) with $P^k_{ij}=\S f$ as
\bq\label{Hess}
\mbox{Hess}(u)(z)(\vec v,\vec v)=S_f(z)(\vec v)\cdot \nabla u(z)+S^0_f(z)(\vec v)u\, ,
\eq
where $S^0_f(z)(\vec v)=\vec v\,^tS^0_f(z)\vec v$.\\

The Bergman metric in the polydisk $\P=\{(z_1,\ldots,z_n): |z_i|<1,\,i=1,\ldots, n\}$ is the hermitian product
defined by the diagonal matrix $$g_{ii}(z)=\frac{2}{(1-|z_i|^2)^2}\; ,$$ see, e.g., \cite{Kr}.
The norm of the Schwarzian operator is defined by
$$
\|S_f(z)\|=\dis\sup_{\|\vec v\|=1}\|S_f(z)(\vec v\,)\|\, ,
$$
where
$$
\|\,\vec v\, \|=\left[\,2\sum_{i=1}^n
\frac{|v_i|^2}{(1-|z_i|^2)^2}\,\right]^{1/2}
$$
\noi
is the Bergman norm of $\vec v \in T_z\P$.
Finally, let
$$
\|S_f\| = \sup_{z\in\P}\|S_f(z)\|\, .
$$
For a system of the form (\ref{system}) that is completely integrable and in canonical form, we define $\|P\|$ to be $\|S_f\|$ of any mapping with Schwarzian
derivatives given by the coefficients $P^k_{ij}$. Let now
$$
\Fa=\{f:\P\rightarrow\Cn\, \mid \; f(0)=0\; , \; Df(0)= \I \; , \; \|S_f\|\leq \al\} \, .
$$

\sm
\noi
It is well known that up to multiplication by a diagonal unitary transformation and a
permutation of the coordinates, the automorphisms of the polydisk are given by
\bq\label{automorphism}
\psi(z)=\psi_a(z)=(\psi_{a_1}(z_1),\ldots,\psi_{a_n}(z_n))\, , \;\;z=(z_1,\ldots,z_n)\in\P\, ,
\eq
where $a=(a_1,\ldots,a_n)\in\P$ and
$\psi_{a_j}(z_j)=\dis\frac{z_j-a_j}{1-\overline{a_j}z_j}\;,\,1\leq
j \leq n\,.$
Because the automorphisms of $\P$ are not M\"{o}bius transformations unless the origin is fixed, the bound on $\|S_f\|$ in $\Fa$ will generally not be preserved when composing with them. This poses a significant difference regarding
the analysis of the corresponding family in $\B$, and it will be useful to estimate the norm of $\|S_{\psi}\|$ for an automorphism of $\P$.

The components $\psi_j=\psi_{a_j}$ of the mapping $\psi=\psi_a$ satisfy
$$
\begin{array}{lll}S^k_{kk}\psi(z) &=&\dis\frac{n-1}{n+1}\,\frac{\psi''_k}{\psi'_k}(z_k)\, ,\\[0.5cm]
S^k_{jj}\psi(z)&=& 0 \; , \;\;  j\neq k\, ,\\[0.3cm]
S^k_{kj}\psi(z)&=& -\dis\frac{1}{n+1}\frac{\psi''_j}{\psi'_j}(z_j)\; ,\;\;  j\neq k\, ,\\
S^k_{ij}\psi(z)&=& 0\; , \; i\neq j\, , \,k\neq j\, .
\end{array}
$$
It follows that for $\vec v=(v_1,\ldots,v_n)$ we have
$$
\vec{v}^{\,t}S^k_{\psi}(z)\vec v=\frac{n-1}{n+1}\cdot \frac{\psi''_k}{\psi'_k}(z_k)v_k^2-\frac{2}{n+1}\sum_{j=1,j\neq k}^n\frac{\psi''_j}{\psi'_j}(z_j)v_kv_j=\frac{v_k}{n+1}\sum_{j=1}^nd_{kj}\, ,
$$
where $d_{kj}=-2\dis\frac{\psi''_j}{\psi'_j}(z_j)v_j$ for $j\neq k$ and $d_{kk}=(n-1)\dis\frac{\psi''_k}{\psi'_k}(z_k)v_k$. Therefore
\bq\label{Ds}
S_\psi(z)(\vec v)=\frac{1}{n+1}\left(v_1D_1,\ldots, v_nD_n\right)\, ,
\eq
where
$$
D_k=\sum_{j=1}^nd_{kj}=(n-1)\frac{\psi''_k}{\psi'_k}(z_k)v_k-2\sum_{j\neq k}\frac{\psi''_j}{\psi'_j}(z_j)v_j \, .
$$
The square of the Bergman norm of $S_\psi(z)(\vec v)$ is given by
$$
A=\frac{1}{(n+1)^2}\sum_{k=1}^n\frac{2|v_k|^2}{(1-|z_k|^2)^2}|D_k|^2 \, ,
$$
and if we assume now that $||\vec{v}||=1$, we see that the sum on the right hand-side is a convex combination of the $|D_k|^2$.
Hence
$$
A\leq \frac{1}{(n+1)^2}\max_i |D_i|^2 \, .
$$
Since $|\psi_i''/\psi_i'(z_i)|\leq 2|a_i|/|1-\overline{a_i}z_i|$ we see from (\ref{automorphism}) that
$$
|D_k|\, \leq \, 2(n-1)\frac{|a_k||v_k|}{|1-\overline{a_k}z_k|}+4\sum_{j\neq k}\frac{|a_j||v_j|}{|1-\overline{a_j}z_j|} \, ,
$$
and using that $\dis \frac{1}{|1-\overline{a_k}z_k|}\leq \frac{2}{1-|z_k|^2}$ together with the Cauchy-Schwarz inequality, we obtain
$$
\begin{array}{lll}
|D_k| & \leq & 4\left(2|a_1|,\dots, (n-1)|a_k|, \ldots, 2|a_n|\right)\cdot\dis\left(\frac{|v_1|}{1-|z_1|^2},\ldots, \frac{|v_n|}{1-|z_n|^2}\right) \\
& & \\
 & \leq & 2\sqrt{2}\sqrt{4|a_1|^2+\cdots+(n-1)^2|a_k|^2+\cdots+ 4|a_n|^2} \\
& & \\
 & \leq & 2\sqrt{2}\sqrt{(n+3)(n-1)}\max_i |a_i| \;\, \leq \;\, 2\sqrt{2}(n+1)\max_i |a_i| \, .
\end{array}
$$

\me
\noi
With  $|z|_{\infty}=\max_i|z_i|$ denoting the polydisk norm, we see that $A\leq 8|\psi(0)|_{\infty}^2$, hence
$\vspace{,1cm}$
$$
\|S_\psi\|\leq 2\sqrt{2}|\psi(0)|_{\infty} \, .
$$

\me
If $f\in \Fa$ and $\psi$ is an automorphism of $\P$, then we see from (\ref{chain rule}) and the invariance of the Bergman metric under $\psi$ that
\bq\label{norm increase}
\|S_{f\circ\psi}\|\leq \|S_f\|+\|S_{\psi}\|\leq \al+2\sqrt{2}|\psi(0)|_{\infty}\, .
\eq

The notation $|z|$ will be kept to refer to the Euclidean norm.
The following result on estimates for the derivative of a holomorphic function in $\D$ will be used in the paper, as well as two lemmas on
specific bounds for the Schwarzians of a mapping $f$ with bounded norm $\|Sf\|$.

\me
\begin{lem}
Let $g$ be holomorphic in $\D$.

\me
\noi
(i) If $\,\dis |g(z)|\leq C$ then $\,\dis |g'(z)|\leq \frac{C}{1-|z|^2}\, .$

\sm
\noi
(ii) If $\,\dis |g(z)|\leq \frac{C}{1-|z|^2}$ then $\,\dis |g'(z)|\leq \frac{4C}{(1-|z|^2)^2}\, . $
\end{lem}

\bp
For the proof, we may assume that $C=1$. The first part then follows at once from the Schwarz-Pick lemma.

\me

For the second part, fix $a\in\D$, which after rotation may be assumed to real and non-negative. The Cauchy integral formula applied on the circle $C_r: \;|\zeta-a|=r<1-a$ gives
$$ |g'(a)|\leq \frac{1}{2\pi}\int_{C_r}\frac{|d\zeta|}{(1-|\zeta|^2)|\zeta-a|^2}
=\frac{1}{2\pi r}\int_0^{2\pi}\frac{d\theta}{1-a^2-r^2-2ar\cos(\theta)} \, .$$
The last integral can be computed using the Weierstrass substitution $\theta=2\arctan(x)$, and gives
$$ |g'(a)|\leq \frac{1}{r\sqrt{U-V}} \, ,$$
where $U=(1-a^2-r^2)^2, \, V=4a^2r^2$. Evaluating the expression $r\sqrt{U-V}$ at the value $r=r(a)$ where it attains its maximum is rather cumbersome, but the simple expression $r=(1-a)/\sqrt{2}$ is close enough. At this value
$$ U-V=\frac{(1-a)^2}{4}\left(a^2+6a+1\right)\, ,$$ and thus
$$
|g'(a)|\leq \frac{2\sqrt{2}}{(1-a)^2\sqrt{a^2+6a+1}}=\frac{2\sqrt{2}(1+a)^2}{(1-a^2)^2\sqrt{a^2+6a+1}} \, .$$
The quantity $(1+a)^2/\sqrt{a^2+6a+1}$ decreases to its minimum value at $a=\sqrt{17}-4$, and is increasing thereafter. From this we can see that the maximum value is attained at $a=1$ and is equal to $\sqrt{2}$. This finishes the proof of the lemma.
\ep

\me
\begin{lem}
If $f\in\Fa$ and $|\vec{v}|_{\infty}=1$ then
\bq\label{Sk bound}
\left|S^k_f(z)(\vec v)\right|\leq \frac{3n\al}{1-|z|_{\infty}^2} \, ,
\eq
\bq\label{S0 bound}
\left|S^0_f(z)(\vec v)\right|\leq \frac{5n^2\al+2n(n+1)\al^2}{(1-|z|_{\infty}^2)^2}\, .
\eq

\end{lem}

\bp For the first part, let $\vec v=(v_1,\ldots,v_n)$ be a vector with $|\vec{v}|_{\infty}=1$. It was shown in \cite{CH1} that when $\|S_f\|\leq\al$ in $\P$, one has $\S f\equiv 0$ for $k\neq i,j$ and $(1-|z_i|^2)|S^k_{ik}f(z)|\leq \sqrt{2}\al$. Hence

\sm
$
\hspace{1cm} \left|S^k_f(z)(\vec v)\right|= |v_k|\dis\left|\left(2\sum_{i=1}^n S^k_{ik}f(z)v_i-S^k_{kk}f(z)v_k\right)\right|
$

\sm
$
\hspace{1.2in} \leq \dis{\frac{2\sqrt{2}\al}{(1-|z_1|^2)}+\cdots \frac{\sqrt{2}\al}{(1-|z_k|^2)}+\cdots \frac{2\sqrt {2}\al}{(1-|z_n|^2)}}
$

\me
$
\hspace{1.2in} =\dis{\frac{2\sqrt{2}(n-1)\al+\sqrt{2}\al}{1-|z|_{\infty}^2}
\,=\, \frac{\sqrt2(2n-1)\al}{1-|z|_{\infty}^2}\,\leq \,\frac{3n\al}{1-|z|_{\infty}^2}}\, .
$

\me
For the second part and the abbreviated notation $\partial_i=\partial/\partial z_i$, it was shown in \cite{He} that
$$
S^0_{ii}f(z)=-\frac{1}{n-1} \sum_{k=1}^n\partial_kS^k_{ii}f(z)+\frac{1}{n-1}\sum_{k,j=1}^nS^k_{ij}f(z)S^j_{ki}f(z)\, ,
$$
and for $i\neq j$
$$
S^0_{ij}f(z)=\partial_jS^i_{ii}f(z)-\partial_iS^i_{ij}f(z)+\sum_{k=1}^nS^k_{ii}f(z)S^i_{kj}f(z)-\sum_{k=1}^nS^k_{ij}f(z)S^i_{ki}f(z)\, .
$$
Since the $\S f$ vanish unless $k=i$ or $k=j$, we find that
$$
S^0_{ii}f(z)=-\frac{1}{n-1}\partial_iS^i_{ii}f(z)+\frac{1}{n-1}\sum_{k=1}^n\left(S^k_{ki}f(z)\right)^2\, .
$$
and for $i\neq j$
$$
S^0_{ij}f(z)=\partial_j S^i_{ii}f(z)-\partial_i S^i_{ji}f(z)-S^j_{ij}f(z)S^i_{ij}f(z)\, .
$$
From Lemma 2.3 we see that
\bq\label{S0ii}
|S^0_{ii}f(z)|\leq \frac{1}{n-1}\frac{4\sqrt{2}\al+2n\al^2}{(1-|z|_{\infty}^2)^2} \, ,
\eq and for $i\neq j$
\bq\label{S0ij}
|S^0_{ij}f(z)|\leq \frac{2\sqrt{2}\al+2\al^2}{(1-|z|_{\infty}^2)^2}\, ,
\eq
because the bounds for $S^i_{ii}$ and $S^i_{ji}$ do not depend on $z_j$ and $z_i$, respectively.
Therefore

\sm
$
\hspace{.7cm} \left|S^0_f(z)(\vec v)\right| = \dis\left|v_1\sum_{k=1}^nS^0_{1k}f(z)v_k
+\cdots+v_n\sum_{k=1}^nS^0_{nk}f(z)v_k\right|
$

\sm
$
\hspace{2.7cm} \leq  \dis n\left(\frac{1}{n-1}\frac{4\sqrt{2}\al+2n\al^2}{(1-|z|_{\infty}^2)^2}+
(n-1)\frac{2\sqrt{2}\al+2\al^2}{(1-|z|_{\infty}^2)^2}\right)
$

\me
$
\hspace{2.7cm} =\dis\frac{2n}{n-1}\frac{\sqrt{2}\left[2+(n-1)^2\right]\al+\left[n+(n-1)^2\right]\al^2}{(1-|z|_{\infty}^2)^2}
\,\leq\, \dis \frac{5n^2\al+2n(n+1)\al^2}{(1-|z|_{\infty}^2)^2}\;\, .
$

\me
\noi
For the last inequality we have used that $2\left[2+(n-1)^2\right]\leq 3n(n-1)$ and $n+(n-1)^2\leq (n+1)(n-1)$ because
$n\geq 2$, and that $3\sqrt{2}\leq5$.

\ep
\begin{lem}
Let (\ref{system}) be a completely integrable system in canonical form with $\|P\|\leq \al$.
For fixed $\zeta_1,\ldots,\zeta_n\in\D$ let $A=(P^k_j)_{k,j}$ be the matrix with components
$$ P^k_j=\zeta_1P^k_{1j}+\cdots+\zeta_nP^k_{nj}\, ,$$
and let $B=(P^0_j)$ be the vector with components
$$P^0_j=\zeta_1P^0_{1j}+\cdots+\zeta_nP^0_{nj}\, .$$
Then with $|A|=\max_{|v|=1}|Av|$ we have
\bq\label{|A|}
(1-|z|_{\infty}^2)|A|\leq 1.6\,n\sqrt{n}\al \, ,
\eq
\bq\label{|B|}
(1-|z|_{\infty}^2)^2|B|\leq c_1(n)\al+c_2(n)\al^2 \, ,
\eq
where
\bq\label{c1,c2}
c_1(n)=\frac{2\sqrt{2n}}{n-1}\left[2+(n-1)^2\right] \;\,, \;\, c_2(n)=\frac{2\sqrt{n}}{n-1}\left[n+(n-1)^2\right] \, .
\eq
Furthermore, if $n\sqrt{n}\al\leq3\sqrt{2}-4=0.2426...$ then
\bq\label{|B| better}
(1-|z|_{\infty}^2)^2|B|\leq 4.5\,n\sqrt{n}\al \, .
\eq
\end{lem}

\bp
We first estimate $|P^k_j|$ for $P^k_{ij}=\S f$. For $k\neq j$ the only term contributing in the sum is $\zeta_kP^k_{kj}$, and thus
$(1-|z|_{\infty}^2)|P^k_j|\leq\sqrt{2}\al$. When $k=j$ all summands may contribute and we see that $(1-|z|_{\infty}^2)|P^k_j|\leq\sqrt{2}n\al$. Therefore, for every column $C$ of $A$ we have
$$
(1-|z|_{\infty}^2)|C|\leq \sqrt{2}\sqrt{n^2+(n-1)}\,\al \, .
$$
The bound (\ref{|A|}) follows from the estimate $\sqrt{2}\sqrt{n^2+n-1}\leq1.6n$.

\sm
For the second part, we use (\ref{S0ii}) and (\ref{S0ij}) to find that
$$
(1-|z|_{\infty}^2)^2|P^0_j|\leq \frac{2}{n-1}\left(\sqrt{2}\left[2+(n-1)^2\right]\al+\left[n+(n-1)^2\right]\al^2\right)\, ,
$$
which establishes (\ref{|B|})-(\ref{c1,c2}) because $|B|\leq\sqrt{n}\max_j|P^0_j|$.

\sm
For the last claim, fix $\epsilon >0$ and suppose that $n\sqrt{n}\al\leq \epsilon$. Then
$$
c_1(n)\al+c_2(n)\al^2\leq \left(c_1(n)+\frac{c_2(n)}{n\sqrt{n}}\epsilon\right)\al \, .
$$
The inequality (\ref{|B| better}) will then hold if
$$
2\sqrt{2}\left[2+(n-1)^2\right]+2\frac{n+(n-1)^2}{n\sqrt{n}}\epsilon \leq 4.5\,n(n-1) \, ,
$$
or equivalently, after dividing by $2$ and shifting the term $\sqrt{2}(n-1)^2$ to the right, if
$$
2\sqrt{2}+\frac{n+(n-1)^2}{n\sqrt{n}}\epsilon \leq(n-1)\left[2.25n-\sqrt{2}(n-1)\right] \, ,
$$
that is, if
$$
\frac{2\sqrt{2}}{n-1}+\frac{n+(n-1)^2}{n(n-1)\sqrt{n}}\epsilon \leq (2.25-\sqrt{2})n+\sqrt{2} \, .
$$
Because the left hand-side is decreasing in $n$ while the right hand-side is increasing, we see that (\ref{|B| better})
will hold provided it holds for $n=2$, that is, if
$$
3\sqrt{2}+\frac34\sqrt{2}\epsilon \leq 4.5 \, ,
$$
which gives the constant $3\sqrt{2}-4$ as the maximum possible value for $\epsilon$. This finishes the proof.
\ep

It will become useful to compare the constant $c_1(n)\al+c_2(n)\al^2$ with the bound $\tau=1.6n\sqrt{n}\al$ in simple terms.
We have that
$$
c_1(n)\al=5\sqrt{2}\,\frac{2+(n-1)^2}{4n(n-1)}\tau\leq\frac{15}{8}\sqrt{2}\tau\leq 3\tau \, ,
$$
while
$$
c_2(n)\al^2=5\,\frac{n+(n-1)^2}{4n(n-1)}\tau\al\leq\frac{15}{8}\tau\al\leq 2\tau\al \, ,
$$
hence
\bq\label{c1+c2}
c_1(n)\al+c_2(n)\al^2\leq(3+2\al)\tau \, .
\eq

\section{Normalization and Finite Norm}

In this section we will establish important bounds for mappings in $\Fa$. In particular, we will show that the family $\Fa$ is compact and will exhibit a connection between the order
$$\laa=\sup_{f\in\Fa}|\nabla J_f(0)|$$
and the size of the largest Euclidean ball centered at the origin that is contained in every image $f(\P)$, $\,f\in\Fa$.
To this end, we will consider the associated family $\Fao$ of mappings $g=T(f)$ for $f\in\Fa$ and $T$ a specific M\"obius transformation that ensures the additional normalization $\nabla J_g(0)=0$. Any such mapping will remain locally biholomorphic away from a possibly existing singular set in $\P$ that depends on $g$ itself. In any case, such a set will stay away some fixed distance $r_0=r_0(n,\al)$ from the origin, and will not exist in $\P$ altogether when $\al$ is sufficiently small. The procedure of normalizing is described in the following lemma.

\begin{lem} Let $f:\P\rightarrow\Cn$ be locally biholomorphic, with $f(0)=0$ and $Df(0)=\I$. Then there exists a M\"obius map $T$ such $g=T\circ f$ satisfies that $g(0)=0$, $Dg(0)=\rm{Id}$ and $\nabla J_g(0)=0.$
\end{lem}

\bp For $a=(a_1,\ldots,a_n)$ let $T$ be a M\"obius map of the form
\bq\label{moebius}
T_a(w)=\left(\frac{w_1}{1+a_1w_1+\cdots+a_nw_n},\ldots,\frac{w_n}{1+a_1w_1+\cdots+a_nw_n}\right)\, .
\eq
The mapping $g=T(f)$ has $g(0)=0$ and $Dg(0)=\I$. Also,  $Dg(z)=DT(f(z))\cdot Df(z)$ hence $J_g(z)=J_T(f(z))\cdot J_f(z)$ and $$\nabla J_g(z)=\nabla J_T(f(z))\cdot Df(z)\cdot J_f(z)+J_T(f(z))\cdot \nabla J_f(z)\,.$$ Evaluating at the origin we obtain that $$\nabla J_g(0)=\nabla J_T(0)+\nabla J_f(0)=-(n+1)(a_1,\ldots,a_n)+\nabla J_f(0)\, ,$$ so that
the lemma holds if $(n+1)(a_1,\ldots, a_n)=\nabla J_f(0)$.
\ep

It is important to observe that  $u=( J_g)^{-\frac{1}{n+1}}=(1+a_1f_1+\cdots+a_nf_n)(J_f)^{-\frac{1}{n+1}}$ is a solution of (\ref{system}) with $P^k_{ij}=\S f$ and the initial conditions $u(0)=1, \nabla u(0)=0$. It vanishes precisely when $1+a_1f_1+\cdots+a_nf_n=0$, that is,  when the image $f(\P)$ meets the hyperplane $1+a_1w_1+\cdots+a_nw_n=0$, or $g$ becomes infinite.

In order to analyze normalization, we will require estimates obtained from comparison techniques for differential equation.

\begin{lem}
Let $h:[0,1)\rightarrow\mathbb{R}$ be a non-negative function satisfying
$$
h''\leq\frac{2a}{1-x^2} h'+\frac{b}{(1-x^2)^2}h\;\, , \; h(0)=1\, , \, h'(0)=0\,,
$$
where $a,b>0$ are constants. Then
\bq\label{h'-a,c}
h'\leq \frac{a+c}{1-x^2}\, h \, ,
\eq
and
\bq\label{h-a,c-1}
h\leq 2\left(\frac{1+x}{1-x}\right)^{\frac{a+c-1}{2}} \, ,
\eq
where $c=\sqrt{1+b+a^2}$.
\end{lem}

\bp
We will analyze the function $h$ as a function of the variable $s=\dis \frac{1}{2}\log\frac{1+x}{1-x}$ defined on
$[0,\infty)$. The change of variable $x=x(s)$ has $x'=1-x^2$. Let
$$
w(s)=\frac{h}{u}(x(s))\, ,
$$
where $u(x)=\sqrt{1-x^2}$. Then $w(0)=1, w'(0)=0$ and $w'=uh'-u'h$, with the agreement that derivatives of $w$ are with respect to $s$, while those of $u,h$ are with respect to $x$. One further derivative yields $w''=(1-x^2)(uh''-hu'')$, and because $u''=-u/(1-x^2)^2$ we see that
$$
w''\leq 2a(uh')+\frac{1+b}{1-x^2}\,uh \, .
$$
We replace the term $uh'$ in the right hand side with $w'+u'h=w'-(x/u)h=w'-xw$, and because $h$, and thus $w$, are non-negative, we obtain
$$
w''\leq 2a(w'-xw)+(1+b)w\leq 2aw'+(1+b)w\, .
$$
Let $\la_{1,2}=-a\pm\sqrt{1+b+a^2}=-a\pm c$ be the roots of the equation $\la(\la+2a)=1+b$, where $\la_1>0>\la_2$. Then
$$
w''+\la_1w'\leq (2a+\la_1)w'+(1+b)w=(2a+\la_1)(w'+\la_1w)=-\la_2(w'+\la_1w) \, ,
$$
which upon amplification by $e^{\la_1s}$ yields
$$
\left(e^{\la_1s}w'\right)'\leq -\la_2\left(e^{\la_1s}w\right)' \, .
$$
The initial conditions $w(0)=1, w'(0)=0$ imply that
$$
e^{\la_1s}w'\leq -\la_2\left(e^{\la_1s}w-1\right) \, ,
$$
hence
\bq\label{w' negative}
e^{\la_1s}(w'+\la_2w)\leq\la_2 \, .
\eq
Multiplication now by $e^{(\la_2-\la_1)s}$ and integration leads to
$$
w\,\leq\, \frac{\la_2e^{-\la_1s}-\la_1e^{-\la_2s}}{\la_2-\la_1}
\,=\,\frac{1}{2c}\left((c-a)e^{(a+c)s}+(c+a)e^{(a-c)s}\right)
$$
$$
=\,\frac{e^{(a+c)s}}{2c}\,\left(c-a+(c+a)e^{-2cs}\right)\, \leq \,e^{(a+c)s} \, .
$$
In terms of the function $h$ we see that
$$
h\leq \sqrt{1-x^2}\left(\frac{1+x}{1-x}\right)^{\frac{a+c}{2}}
=(1+x)\left(\frac{1+x}{1-x}\right)^{\frac{a+c-1}{2}} \leq 2\left(\frac{1+x}{1-x}\right)^{\frac{a+c-1}{2}} \, .
$$
Finally, we see from (\ref{w' negative}) that $w'+\la_2w\leq 0$, which yields
$$
h'u+\frac{x}{u}h\leq (a+c)\frac{h}{u} \, ,
$$
hence $$
h'\leq \frac{a+c}{1-x^2}\,h \, .
$$
This proves (\ref{h'-a,c}) and (\ref{h-a,c-1}).
\ep

\begin{lem}
Let (\ref{system}) be a completely integrable system in canonical form with $\|P\|\leq \al$. Then the solution $u$ with $u(0)=1, \nabla u(0)=0$ satisfies
\bq\label{|u| bound}
|u|\leq 2\left(\frac{1+|z|_{\infty}}{1-|z|_{\infty}}\right)^{\ga}\, ,
\eq
\bq\label{|Du| bound}
|\nabla u|\leq 2\,\frac{1+2\ga}{1-|z|_{\infty}^2}
\left(\frac{1+|z|_{\infty}}{1-|z|_{\infty}}\right)^{\ga}\, ,
\eq
where $\dis 2\ga=\sqrt{1+b+a^2}-1$, $a=0.8n\sqrt{n}\al$ and $b=2(3+2\al)\sqrt{n}a$.
\end{lem}

\bp
For $\zeta\in\partial\P$ we let $v(t)=u(t\zeta)$, $t\in[0,1]$. We seek to estimate the functions
$w_j(t)=\partial u/\partial z_j(t\zeta)$. The $t$-derivatives of these quantities are given by
$$
\begin{array}{lll}
w_j'&=& \dis\frac{\partial^2u}{\partial z_1\partial z_j}\zeta_1+\cdots + \frac{\partial^2u}{\partial z_n\partial z_j}\zeta_n\\[0.5 cm] &=& \dis \left(\sum_{k=1}^n P^k_{1j}\frac{\partial u}{\partial z_k}+P^0_{1j} u\right)\zeta_1
+\cdots+\dis \left(\sum_{k=1}^n P^k_{nj}\frac{\partial u}{\partial z_k}+P^0_{nj}u\right)\zeta_n\\[0.5cm]
&=&\dis \sum_{k=1}^n\left(P^k_{1j}\zeta_1+\cdots P^k_{nj}\zeta_n\right) \frac{\partial u}{\partial z_k}
+\left(P^0_{1j}\zeta_1+\cdots P^0_{nj}\zeta_n\right)u\\[0.5cm]
&=&\dis \sum_{k=1}^n\left(P^k_{1j}\zeta_1+\cdots P^k_{nj}\zeta_n\right) w_k
+\left(P^0_{1j}\zeta_1+\cdots P^0_{nj}\zeta_n\right)u=\sum_{k=1}^nP^k_jw_k+P^0_ju \, .
\end{array}
$$
Let $\varphi(t)$ the vector function in $\C^{n}$ with coordinates $\varphi_i=w_i$ for $i=1,\ldots, n$.
Since $v'=\nabla u\cdot\zeta$ we can see that
$$
\varphi'=A\cdot \varphi+uB\, ,
$$
where $A=(P^k_j)_{k,j}$ is an $n\times n$ matrix and $B$ is the vector with components $P^0_j$.
Since $|\varphi|'\leq |\varphi'|$ we see that
$$
|\varphi|'\leq |A||\varphi|+|B||u| \, ,
$$
where $|A|=\max_{|v|=1}|Av|$. The function $\phi(t)=\dis\sqrt{n}\int_0^t|\varphi(s)|ds +1$ is an upper bound for $|u|$ with $\phi(0)=1, \phi'(0)=0$
and satisfies
$$
\phi''\leq|A|\phi'+\sqrt{n}|B|\phi  \, .
$$
\me
\noi
The inequalities (\ref{|A|}) and (\ref{c1+c2}) show that $(1-t^2)|A|\leq 2a$ and
$(1-t^2)^2|B|\leq b=2\sqrt{n}(3+2\al)a$, where $a=0.8n\sqrt{n}\al$. It follows now from Lemma 3.2 that
$$
\hspace{.8cm} \phi\leq 2\left(\frac{1+t}{1-t}\right)^{\ga}\;\; , \quad\;
\phi'\leq \,2\, \frac{1+2\ga}{1-t^2}\left(\frac{1+t}{1-t}\right)^{\ga}\, ,
$$
for $2\ga=a+\sqrt{1+b+a^2}-1$. This finishes the proof.

\ep

\begin{cor}
Let $\al>0$ be fixed. The set of solutions $u$ with $u(0)=1, \nabla u(0)=0$ of completely integrable systems (\ref{system}) in canonical form with $\|P\|\leq \al$ is compact.
\end{cor}

\begin{lem} Let $g\in \Fao$. Then there exists $r_0=r_0(n,\al)>0$ such that $g$ is holomorphic in $\B(0,r_0)$.
\end{lem}

\bp Let $g=T(f)\in \Fao$ with $f\in \Fa$ as before.
Then  $u=( J_g)^{-\frac{1}{n+1}}=(1+a_1w_1+\cdots+a_nw_n)(J_f)^{-\frac{1}{n+1}}$ is a solution of (\ref{system}) with $P^k_{ij}=\S f, P^0_{ij}=S^0_{ij}$ that vanishes precisely when $g$ ceases to be holomorphic. We will show that $u\neq 0$ in some neighborhood of the origin.

\sm
For $\zeta\in\partial\P$ we define $v(t)=|u(t\zeta)|$, $t\in[0,1)$. Therefore
$vv'=\mbox{Re}\{ \overline{u}(\nabla u\cdot \zeta)\}$ and
$$
vv''+(v')^2=|\nabla u\cdot\zeta|^2+\mbox{Re}\{\overline{u}\,\mbox{Hess}(u)(\zeta,\zeta)\} \, .
$$
Since $|\nabla u\cdot\zeta|^2-(v')^2\geq 0$ we conclude from (\ref{Hess}) that
$$
v(t)v''(t)\geq \mbox{Re}\left\{\overline{u}S_f(t\zeta)(\zeta)\cdot\nabla u(t\zeta)+|u|^2S^0_f(t\zeta)(\zeta)\right\}\, ,
$$
and hence
\bq\label{v'' bound}
v''(t) \geq -|S_f(t\zeta)(\zeta)||\nabla u(t\zeta)|-|S^0_f(t\zeta)(\zeta)|v(t)\, . \hspace{1cm}
\eq
The bound (\ref{|Du| bound}) and Lemma 2.4 give that
$$
v''+\frac{\varepsilon}{(1-t^2)^2} v\geq-\frac{\delta}{(1-t^2)^2}\left(\frac{1+t}{1-t}\right)^{\ga}\, ,
$$
with $v(0)=1, v'(0)=0$, and $\varepsilon=\varepsilon(n,\al), \delta=\delta(n,\al)$ are positive constants. We rewrite this inequality as
$$
v''+\frac{\varepsilon}{(1-t^2)^2} v=-\frac{\delta}{(1-t^2)^2}\left(\frac{1+t}{1-t}\right)^{\ga}+\eta(t) \, ,
$$
for some $\eta(t)\geq 0$, and claim that $v\geq y$ for $t\in[0,t_0]$, where $t_0$ is the first zero of the solution $y=y(t)$ of
$$
y''+\frac{\varepsilon}{(1-t^2)^2} y=-\frac{\delta}{1-t^2}\left(\frac{1+t}{1-t}\right)^{\ga},\quad y(0)=1\,,\; y'(0)=0\, .
$$
For this, observe first that $y''+\varepsilon(1-t^2)^{-2}y\leq 0$, hence from Sturm comparison we have that $0\leq y\leq h$ on $[0,t_0]$ for the solution of
$$
h''+\frac{\varepsilon}{(1-t^2)^2}h= 0\; \,, \;h(0)=1 \, ,\, h'(0)=0 \, .
$$
For $s>0$ we consider the solution
$$
v_s''+\frac{\varepsilon}{(1-t^2)^2} v_s=-\frac{\delta}{(1-t^2)^2}\left(\frac{1+t}{1-t}\right)^{\ga}+\eta(t) \;, \;v_s(0)=1+s \, ,\, v_s'(0)=0 \, .
$$
The functions $\{v_s\}_{s>0}$ will converge to $v$ locally uniformly on $[0,1)$ as $s\rightarrow 0$. Since difference $w=v_s-y$ satisfies
$$
w''+\frac{\varepsilon}{(1-t^2)^{2}}w=\eta(t)\geq 0 \; \, , \;w(0)=s \, ,\, w'(0)=0 \, ,
$$
we see, once more, from Sturm comparison that $w\geq sh\geq 0$ on $[0,t_0]$. Therefore, $v_s\geq y$ on $[0,t_0]$, which proves the claim by letting $s\rightarrow 0$.
We conclude finally that $u(z)\neq 0$ for all $z\in\P$ with $|z|<r_0=t_0$, which completes the proof.
\ep

\begin{thm}
The family $\Fa$ has finite order bounded by $1/s_0$ for some $s_0=s_0(n,\al)>0$. The image of every mapping in $\Fao$ covers a Euclidean ball of radius $s_0$ centered at the origin, and every mapping in $\Fa$ covers a corresponding Euclidean ball of radius $s_0/2$ centered at the origin.

\end{thm}

\bp
Corollary 3.4 and Lemma 3.5 show that the family $\Fao$ is compact in any ball $B(0,t_1)$ for $t_1<t_0$. We conclude that there exists $s_0>0$ such that $B(0,s_0)\subset g(\P)$ for every $g\in\Fao$. An estimate for $s_0$ in terms of $n, \al$ will be given in the final section. The mapping
$$
f=T^{-1}(g)=\left(\frac{g_1}{L(g)},\cdots,\frac{g_n}{L(g)}\right)\;\, , \; L(g)=1-(a_1g_1+\cdots +a_ng_n)
$$
is holomorphic in $\P$, hence $L(g)=1-(a_1g_1+\cdots +a_ng_n)$ cannot vanish in $\P$. Since $g$ covers a Euclidean ball of radius $s_0$ it follows that
$$
\frac{1}{(n+1)^2}|\nabla J_f(0)|^2=|a_1|^2+\cdots+|a_n|^2\leq \frac{1}{s_0^2} \, .
$$
This implies that
$$
\sup_{f\in\Fa}|\nabla J_f(0)|\leq \frac{n+1}{s_0}<\infty \, .
$$
Finally, let $f\in \Fa$ and let $g=T(f)$ be the corresponding mapping in $\Fao$ as before. For $z\in U=g^{-1}(B(0,s_0))$ we see from the triangle and the Cauchy-Schwarz inequalities that
$$|L(g)(z)|=|1-a_1g_1-\cdots-a_ng_n|\leq 1+|a_1g_1|+\cdots+|a_ng_n|\leq 2 \, .$$
It follows that
$$
|f(z)|=\frac{|g(z)|}{|L(g)|}\geq \frac{|g(z)|}{2} \, ,
$$
therefore
$$
\liminf_{z\rightarrow\partial U}|f(z)|\geq \frac{s_0}{2} \, .
$$
\ep

\sm
As it turns out, for $\al\leq\al_0(n)$ small enough the mappings in $\Fao$ will be holomorphic in all of $\P$. Equivalently,
the solutions $u$ of (\ref{system}) with $u(0)=1, \nabla u(0)=0$ do not vanish in $\P$. To establish this we will require
a variation of Lemma 3.3 focusing on small $\al$. The phenomenon of normalizing at the origin while keeping the mapping regular
also occurs in one variable, with $\mathcal{F}_{\al_0}$ corresponding to the Nehari class $||Sf||\leq 2$. Indeed, for a Nehari mapping $f=z+a_2z^2+\cdots$ it is known that
$-1/a_2\notin f(\D)$, and thus the mapping
$$
g=\frac{f}{1+a_2f}
$$
which is normalized to have $g''(0)=0$ remains holomorphic in $\D$. It is to be noted that the mappings in $\Fa$ are univalent if $3\sqrt{2}\al\leq 1$ \cite{CH1}.

\sm

\begin{lem}
Let $h$ be the solution of the equation
\bq\label{1.6-4.5}
h'=\frac{1.6c}{1-x^2}h+\frac{4.5c}{(1-x^2)^2}+h^2\;\, , \; h(0)=0\, ,
\eq
where $c>0$ is a constant. If $c\leq \dis\frac{1}{6.1}=0.163...$ then the solution $h$ exists on the entire interval $[0,1)$ and
\bq\label{h-y}
h\leq \frac{x}{1-x^2} \, .
\eq
The bound for $c$ cannot be improved beyond $\dis\frac{2}{6.1+\sqrt{(6.1)^2-(1.6)^2}}=0.166...$.
\end{lem}

\bp
Because of the term $h^2$ in the differential equation, $h$ will become infinite before reaching $x=1$ if at any point
$x_0<1$ the value $h(x_0)$ is too large. This will happen for large $c$ because of the term $4.5c/(1-x^2)^2$ on the right-hand side.

The function $y=\dis \frac{x}{1-x^2}$ satisfies
\bq\label{y'}
y'=\frac{1+x^2}{(1-x^2)^2}\geq\frac{1.6c}{1-x^2}y+\frac{4.5c}{(1-x^2)^2}+y^2
\eq
provided that
$$
1+x^2\geq 1.6cx+4.5c+x^2 \, .
$$
This inequality will hold for all $x\in[0,1)$ precisely when $6.1c\leq 1$. It remains to be shown that $h\leq y$.
To that end, we rewrite (3.10) as
$$
y'=\frac{1.6c}{1-x^2}y+\frac{4.5c}{(1-x^2)^2}+y^2+\eta(x) \, ,
$$
for some $\eta(x)\geq 0$ on $[0,1)$, and we consider the solution $w_s$ of this last equation with initial condition $w_s(0)=s>0$. The function $w_s$ will exist on some interval $[0,x_s)$ for some $x_s<1$ (or even $x_s$=1), but will converge to $y$ locally uniformly on $[0,1)$ as $s\rightarrow 0$. Hence $x_s$ will converge to $1$. On $[0,x_s)$ we have that
$$
(w_s-h)'=\frac{1.6c}{1-x^2}(w_s-h)+(w_s^2-h^2)+\eta(x)
$$
$$
\geq \frac{1.6c}{1-x^2}(w_s-h)+(w_s+h)(w_s-h) \, .
$$
Since evidently $w_s, h\geq 0$, we conclude that $w_s-h>0$ because $(w_s-h)(0)=s>0$. The inequality $y\geq h$ results from letting $s$ tend to $0$.

\me
For the last claim, we express (\ref{1.6-4.5})in terms of the function $\phi=(1-x)h$, as
$$
(1-x)\phi'=-\left(1-\frac{1.6c}{1+x}\right)\phi+\frac{4.5c}{(1+x)^2}+\phi^2 =(\phi-\beta)^2+\gamma \, ,
$$
where $\dis\beta=\frac12\left(1-\frac{1.6c}{1+x}\right)$ and
$$
\ga=\frac{4.5c}{(1+x)^2}-\beta^2\, .
$$
The quantity $\gamma$ will be decreasing in $x$ for, say, $c\leq 1$, and thus
$$\ga\geq \frac{4.5c}{4}-\frac14(1-0.8c)^2:=d \, .
$$
The number $d$ will be positive provided $0.64c^2-6.1c+1<0$, which happens for $c\in(0.166...,9.364...)$.
Thus for any $c\in(0.166...,1)$ we have $d>0$, so that $\phi$ is strictly increasing. Hence $\phi'\geq d(1-x)^{-1}$ and there exists $x_0\in(0,1)$ so that $\phi(x_0)>1/2>\beta$. For $x\geq x_0$ we thus have that $\phi>1/2>\beta$ and
$$
\frac{\phi'}{(\phi-\frac12)^2}\geq \frac{1}{1-x} \, ,
$$
which upon integration shows that $\phi$ becomes infinite at some $x_1<1$.
\ep

\begin{thm}
Let (\ref{system}) be a completely integrable system in canonical form with $\|P\|\leq \al$, and let $u$ be the solution with $u(0)=1, \nabla u (0)=0$.  If
\bq\label{small al}
n\sqrt{n}\al \leq \frac{1}{6.1}
\eq
then
$u$ does not vanish in $\P$ and
\bq\label{|Du/u|}
\left|\frac{\nabla u}{u}\right|\leq\frac{|z|_{\infty}}{1-|z|_{\infty}^2} \, .
\eq
\end{thm}

\bp
Let $u$ be as stated in the theorem. For fixed $\zeta\in\partial\P$ we let
$$
\va(t)= \frac{\nabla u}{u}(t\zeta)\, ,
$$
which is defined on some maximal interval $[0,t_0)$. We claim that $t_0=1$ when $\al$ is sufficiently small.
We see that
$$
\va'(t)=\left[\frac1u\mbox{Hess }u(\zeta,\cdot)\right](t\zeta)-\left[\frac{\nabla u\cdot\zeta}{u^2}\nabla u\right](t\zeta)\, ,
$$
It follows from (\ref{system}) that
$$
\mbox{Hess }u(\zeta,\cdot) =A\cdot\nabla u+uB\, ,
$$
where $A=(a_{ij})$ is the $n\times n$ matrix given by $a_{ij}=\sum_{k=1}^nS_{ik}^jf\zeta_k$, and $B=(b_i)$ is the vector with
$b_i=\sum_{k=1}^nS^0_{ik}\zeta_k$. Hence the function $g(t)=|\va(t)|$ has $g(0)=0$ and satisfies
\bq\label{g'}
g'\leq|\va'|\leq \frac{ag}{1-t^2}+\frac{b}{(1-t^2)^2}+g^2 \, ,
\eq
where, as before, we use (2.11) and (2.14) that see that $a\leq 1.6n\sqrt{n}\al$ and $b\leq 4.5n\sqrt{n}\al$. The theorem now follows from Lemma 3.7 with $c=n\sqrt{n}\al$.
\ep

We draw the following immediate corollary.

\begin{cor}
If $n\sqrt{n}\al\leq \dis\frac{1}{6.1}$ then all mappings in $\Fao$ are holomorphic in $\P$.
\end{cor}

The next result is also a direct consequence of Theorem 3.8.

\begin{thm}
Let (\ref{system}) be a completely integrable system in canonical form with $\|P\|\leq \al$, and let $u$ be the solution with $u(0)=1, \nabla u (0)=0$. If
$$
n\sqrt{n}\al \leq \frac{1}{6.1}
$$
then
$$
\left(1-|z|_{\infty}^2\right)^{\frac{\sqrt{n}}{2}}\leq |u|\leq \frac{1}{\left(1-|z|_{\infty}^2\right)^{\frac{\sqrt{n}}{2}}}  \, .
$$
\end{thm}

\bp
Since $u(0)=1$, we find from integrating the bound in (\ref{|Du/u|}) that
$$
|\log u|\leq \frac{\sqrt{n}}{2}\,\log\frac{1}{1-|z|_{\infty}^2}\, ,
$$
which gives the result.
\ep

The following corollary follows at once from the previous theorem, and mimics results for mappings in linearly invariant families \cite{GK}.

\begin{cor}
Let $f\in \Fao$ for $n\sqrt{n}\al\leq \dis\frac{1}{6.1}$. Then
$$
(1-|z|_{\infty}^2)^{\frac{\sqrt{n}(n+1)}{2}}\leq |J_f|\leq \frac{1}{(1-|z|_{\infty}^2)^{\frac{\sqrt{n}(n+1)}{2}}}  \, .
$$
\end{cor}

\begin{thm}\label{cota-u_0} Let $f\in \Fa$. For $z\in \P$ we let $\eta(z)=\min_{|v|=1}|Df(z)v|$ and $\beta=\al+2\sqrt{2}|z|_{\infty}$.
Then the image $f(\P)$ covers a Euclidean ball centered at $f(z)$ of radius
at least
$$
\frac12(1-|z|_{\infty})\eta(z)s_0(n,\beta) \, ,
$$
and
$$
|\nabla J_f(z)|_{\infty}\leq \frac{\la_{\beta}+2|z|_{\infty}}{1-|z|_{\infty}^2}\left|J_f(z)\right| \, ,
$$
where $s_0(n,\beta)$ is the constant defined in Theorem 3.6.
\end{thm}

\bp

For $z_0\in \P$ fixed, let $g$ be defined by
$$
g(z)=[D\psi(0)]^{-1}[Df(z_0)]^{-1}(f(\psi(z))-f(z_0)),\quad \psi(0)=z_0\, ,
$$
where $\psi=\psi_{z_0}$ is an automorphism as given by (\ref{automorphism}). Then $g(0)=0$, $Dg(0)=\I$ and we have that
$\|S_g\|\leq \al+2\sqrt{2}|z_0|_{\infty}$. Hence by (\ref{norm increase}) $g\in\mathcal F_{\beta}$ for $\beta=\al+2\sqrt{2}|z_0|_{\infty}$, and by Theorem 3.6, $g(\P)$ contains
a ball centered at the origin of Euclidean radius at least $\frac12 s_0(n,\beta)$. Since $D\psi(0)$ is a diagonal matrix of elements $1-|z^0_k|^2$,
we see that $\min_{|v|=1}|D\psi(0)v|\geq 1-|z_0|_{\infty}^2$. The first statement of the theorem follows now from the definition of $g$.
\sm
On the other hand,
$$
J_g(z)=J_\psi(0)^{-1}\cdot J_f(z_0)^{-1}\cdot J_f(\psi(z))\cdot J_\psi(z)\, ,
$$
which implies that $\log J_g(z)=\log J_f(\psi(z))+ \log J_\psi(z)+\mbox{constant}$, and
\bq\label{jacobian}
\nabla J_g(0)=\frac{\nabla J_f}{J_f}(z_0)D\psi(0)+\frac{\nabla J_\psi}{J_\psi}(0)\, .
\eq
Because $g\in\mathcal F_{\beta}$ we must have that $|\nabla J_g(0)|\leq \la_{\beta}$.  In addition, with $z_0=(a_1,\ldots,a_n)$
$$
J_\psi(z)=\prod_{j=1}^n\frac{1-|a_j|^2}{(1-\overline{a_j}z_j)^2}\, ,
$$
therefore
$$
\frac{\nabla J_\psi}{J_\psi}(z)=\left(\frac{2\overline{a_1}}{1-\overline{a_1}z_1},
\ldots,\frac{2\overline a_n}{1-\overline{a_n}z_n}\right) \, ,
$$
so that $\dis \frac{\nabla J_\psi}{J_\psi}(0)=2\overline z_0$.
From (\ref{jacobian}) we obtain
$$
\left|\frac{\nabla J_f}{J_f}(z_0)D\psi(0)\right|_{\infty}\leq \lambda_{\beta}+2|z_0|_{\infty}\, .
$$
Since $D\psi(0)$ is a diagonal matrix with elements $1-|a_k|^2$, we have
$$
(1-|z_0|_{\infty}^2)|\nabla J_f(z_0)|_{\infty}\leq (\lambda_{\beta}+2|z_0|_{\infty})|J_f(z_0)|\, ,
$$
which shows that
$$
|\nabla J_f(z_0)|_{\infty}\leq \frac{\lambda_{\beta}+2|z_0|_{\infty}}{1-|z_0|_{\infty}^2}\left|J_f(z_0)\right|\, .
$$
\ep

\section{Order}

The purpose in this section is to establish an explicit upper bound for the order of the family $\Fa$. We will consider the families $\Far$ consisting of all mappings $f$ that are defined and locally biholomorphic for $|z|_{\infty}<r<1$, which have $f(0)=0, Df(0)=\I$, and $\|\S f\|\leq \al$; here $\|\cdot\|$ stands for the Bergman norm in $\P$ restricted to the set $|z|_{\infty}<r$. We see that $\Fa\subset\Far$, and it follows as before that $\Far$ is compact. For fixed $r$ we see that
\bq\label{mur}
\mar=\sup_{f\in\Far}{|\nabla J_f(0)|} <\infty \, ,
\eq
which is attained for some mapping in the class. In addition, $\la_{\al}\leq\mu_r(\al)$.

\begin{thm}
Suppose $r^2<1/5$. Then the function $\mu_r=\mar$ is locally Lipschitz for $\al>0$, and whenever the derivative exists,
\bq\label{lipschitz}
\al\mu_r'(\al)\leq C(r)\mu_r(\al)\, ,
\eq
with $\dis C(r)=\frac{1-r^2}{1-5r^2}$.
\end{thm}

\bp
Let $\al>0$ be fixed, and let $f\in\Far$ be extremal for (\ref{mur}). We consider the dilation of $f$ given by $ g(z)=\frac1sf(sz)$ for $0<s<1$. It is easy to see that $|\nabla J_g(0)|=s|\nabla J_f(0)|=s\mu_r(\al)$, and we claim that $g\in \mathcal{F}_{\beta,r}$ for
$$\beta=s\left(\frac{1-s^2r^2}{1-r^2}\right)^2\al <\al \, .$$
From the composition law (\ref{chain rule}) we see that $\S g(z)=s\S f(sz)$, hence in order to estimate $\|S_g(z)\|$ we must consider how the Bergman norm of given vector $\vec{v}$ changes when measured at $z$ and $sz$. Let $\vec{v}=(v_1,\ldots,v_n)$ be a Bergman unit vector at $z=(z_1,\ldots,z_n)$. The square of the Bergman norm $\|\vec{v}\|_{sz}$ at $sz$ is given by
$$
\sum_i\frac{2|v_i|^2}{(1-s^2|z_i|^2)}=\sum_it_i\left(\frac{1-|z_i|^2}{1-s^2|z_i|^2}\right)^2 \, ,
$$
where $t_i=2|v_i|^2/(1-|z_i|^2)^2$ are non-negative numbers with $\sum_it_i=1$. It follows that
$\|\vec{v}\|^2_{sz}$ is a number between the maximum and minimum value of the quantities $\left[(1-|z_i|^2)/(1-s^2|z_i|^2)\right]^2$.
Because the function $(1-x^2)/(1-s^2x^2)$ is decreasing for $x\in[0,1)$ we see that
$$
\left(\frac{1-r^2}{1-s^2r^2}\right)^2 \leq \|\vec{v}\|_{sz}^2 \leq 1 \, ,
$$
therefore
$$
\|S_g(z)\|\leq s\left(\frac{1-s^2r^2}{1-r^2}\right)^2\al \, .
$$
Simple algebra shows that $s\left(\frac{1-s^2r^2}{1-r^2}\right)^2<1$ for $r^2<1/5$, hence $g\in \mathcal{F}_{\beta,r}$ and $\beta<\al$ as claimed.

Finally, since $|\nabla J_g(0)|=s\mar$ we conclude that $\mu_r(\beta)\geq s\mar$, and so
$$\mar-\mu_r(\beta)\leq (1-s)\mar\leq C(r)(\al-\beta)\frac{\mar}{\al} \, .$$
This proves the theorem.
\ep

For our next theorem, recall the order $\laa=\sup_{f\in\Fa}|\nabla J_f(0)|$.

\begin{thm}
For each $\si>0$ there exists $m(\si)>0$ such that
$$
\laa\leq 2(n+1)+m(\si,n)\al^{1+\si} \, ,
$$
for some $m(\si,n)$ depending on $\si$ and $n$.
\end{thm}

\bp
Let $\al_1>0$ be fixed. It follows from (\ref{lipschitz}) that for $\al\geq \al_1$
$$
\mar\leq \mu_r(\al_1)\left(\frac{\al}{\al_1}\right)^{C(r)} \, ,
$$
where $C(r)=(1-r^2)/(1-5r^2)>1$ approaches $1$ as $r\rightarrow 0$. Hence, for $\al\geq \al_1$ we also have
\bq\label{growth}
\laa \leq \frac{\mu_r(\al_1)}{\al_1^{C(r)}}\,\al^{C(r)} \, .
\eq
We now estimate $\la_0$, that is, the order of the family $\mathcal{F}_0$ consisting of M\"{o}bius mappings of the form
$$
f(z)=\left(\frac{z_1}{l(z)}, \ldots, \frac{z_n}{l(z)}\right) \; , \;\, l(z)=1-a_1z_1-\cdots-a_nz_n \, ,
$$
with the requirement that $l(z)\neq 0$ in $\P$. This happens iff $\sum_i|a_i|\leq 1$. Hence
$$
|\nabla J_f(0)|=(n+1)\sqrt{|a_1|^2+\cdots+|a_n|^2}\leq(n+1)\left(\sum_i|a_i|\right)^2\leq n+1 \, ,
$$
with equality when a single $a_i$ has absolute value 1. This shows that $\la_0=n+1$,
and the theorem follows now from (\ref{growth}) by choosing $\al_1=\al_1(n)$ so that, for example, $\la_{\al_1}=2\la_0$.

\ep

\bi
\noi
{\bf Acknowledgement:} We are grateful to the referee for a very thorough revision of the original manuscript and the first revision. The detailed reports have allowed for several improvements of the paper.

\bi
\noi
{\small Facultad de Matem\'aticas, Pontificia Universidad Cat\'olica de Chile,
Casilla 306, Santiago 22, Chile,\, \email{mchuaqui@mat.uc.cl}

\me
\noi
Facultad de Ingenier\'ia y Ciencias,
Universidad Adolfo Ib\'a\~nez,
Av. Padre Hurtado 750, Vi\~na del Mar, Chile,\,
\email{rodrigo.hernandez@uai.cl}

\end{document}